\documentclass{amsart}

\usepackage{amsmath}
\usepackage{amsfonts}
\usepackage{amssymb} 
\usepackage{amsthm} 
\usepackage{mathptmx}
\usepackage{latexsym}

\numberwithin{equation}{section}

\scrollmode

\newtheorem{thm}{Theorem}[section]
\newtheorem{lem}[thm]{Lemma}
\newtheorem{prop}[thm]{Proposition}

\theoremstyle{definition}
\newtheorem{defn}[thm]{Definition}

\theoremstyle{remark}
\newtheorem{rem}[thm]{Remark}

\newcommand{\al}{Amitsur-Levitzki }

\newcommand{\osp}{\mathfrak{osp}}
\newcommand{\spk}{\mathfrak{sp}}
\newcommand{\gl}{\mathfrak{gl}}
\newcommand{\ok}{\mathfrak{o}}
\newcommand{\Sk}{\mathfrak{S}}
\newcommand{\g}{\mathfrak{g}}
\newcommand{\gt}{\tilde{\g}}
\newcommand{\hk}{\mathfrak{h}}
\newcommand{\lk}{\mathfrak{l}}

\newcommand{\Xc}{\mathcal{X}}
\newcommand\Fc{\mathcal F}
\newcommand\Lc{\mathcal L}

\newcommand\Ac{\mathcal A}

\newcommand\Pc{\mathcal P}

\newcommand{\Ab}{\mathsf{A}}

\newcommand\CC{\mathbb C}

\newcommand\ZZ{\mathbb Z}

\newcommand{\Tr}{\operatorname{Tr}}
\newcommand{\ze}{{\scriptscriptstyle{\overline{0}}}}
\newcommand{\un}{{\scriptscriptstyle{\overline{1}}}}
\newcommand\gO{\g_{\ze}}
\newcommand\gI{\g_{\un}}
\newcommand\VO{V_{\ze}}
\newcommand\VI{V_{\un}}

\newcommand{\So}{\operatorname{S}}
\newcommand{\Ao}{\operatorname{A}}
\newcommand{\Do}{\operatorname{D}}
\newcommand{\Lo}{\operatorname{L}}
\newcommand{\Ro}{\operatorname{R}}
\newcommand{\io}{\operatorname{\iota}}

\newcommand{\ad}{\operatorname{ad}}

\newcommand\dsum{\displaystyle \sum}
\renewcommand\dfrac{\displaystyle \frac}
\newcommand\ds\displaystyle
\newcommand\dc[1]{[\![#1]\!]}

\newcommand\e{\varepsilon}

\renewcommand\ge{\geqslant}
\newcommand{\End}{\operatorname{End}}
\newcommand{\ext}{\operatorname{Ext}}
\newcommand{\sym}{\operatorname{Sym}}
\newcommand{\Id}{\operatorname{Id}}
\renewcommand\hat\widehat
\renewcommand\tilde\widetilde
\newcommand{\spa}{\operatorname{span}}
\newcommand{\str}{\operatorname{str}}
\newcommand\stimes{\mathop{\otimes}\limits_s}
\newcommand\adot{\mathop{\cdot}\limits_a }
\newcommand\sdot{\mathop{\cdot}\limits_s }

\begin{document}

\title[Back to the Amitsur-Levitzki theorem]{Back to the Amitsur-Levitzki theorem: a super version for the orthosymplectic Lie superalgebra $\osp (1, 2n)$} 

\author{Pierre-Alexandre Gi\'e, Georges Pinczon, Rosane Ushirobira}

\address{Institut de Math\'ematiques de Bourgogne, Universit\'e de
  Bourgogne, B.P. 47870, F-21078 Dijon Cedex, France} 

\email{pagie, gpinczon, rosane@u-bourgogne.fr}



\begin{abstract} 
We prove an \al type theorem for the Lie superalgebras \linebreak $\osp(1,2n)$
inspired by Kostant's cohomological interpretation of the classical
theorem. We show that the Lie superalgebras $\gl(p,q)$ cannot satisfy an \al
type super identity if $pq \neq 0$ and conjecture that neither can any other
classical simple Lie superalgebra with the exception of $\osp(1,2n)$.
\end{abstract}

\maketitle


\section{Introduction}

The \al theorem states that $\gl (n)$ satisfies the standard polynomial
identity of order $2n$. More precisely:

\medskip

{\bf{\em{Theorem:}}} Define for $X_i \in \gl (n), k \geq 1$: 

\[ I_k(X_1, \dots, X_k) := \sum_{\sigma \in \Sk_k} \epsilon (\sigma)
X_{\sigma(1)} \dots X_{\sigma(k)}. \] \hskip6mm Then $I_{2n} = 0$.
\footnote{It is easy to see that $I_k \neq 0$ if $k < 2n$ and from $I_{2n} =
  0$, that $I_k = 0$ if $k > 2n$ \cite{Jacobson}.}

\medskip

\noindent Amitsur and Levitzki proved their theorem using an inductive method
that does not explain why such identity exists \cite{Amitsur-Levitzki,
Jacobson}. Later, several simplifications and improvements of their proof,
including graphical ones and several new proofs, were given \cite{Kostant58,
Swan, Jacobson, Rowen74, Rosset, Rowen80, Kostant81}. However, all of these
proofs but Kostant's lack a real interpretation of the result.

Eight years after \al, B. Kostant published a truly beautiful proof of their
theorem, based on the cohomology of Lie algebras \cite{Kostant58}. Besides
explaining the existence of the theorem, Kostant proved with his method that
$\ok (2n)$ satisfies the standard polynomial identity of order $4n -2$ (as a
consequence of the particular structure of its invariants due to the existence
of the Pffafian). Another proof of this result was later obtained by Rowen
using a direct method, but with some difficulties \cite{Rowen80}. Finally in
1981 \cite{Kostant81}, Kostant closed the subject once and for all by
providing a very nice interpretation of the theorem in the context of
representation theory and generalizing it using his separation of variables
theorem \cite{Kostant63}. To our knowledge, no one has returned to the \al
theorem since then.

A few comments can be made about Kostant's proofs of the \al theorem. First,
both proofs use the polynomial structure of the ring of invariants of a semi
simple Lie algebra. Second, his cohomological proof is based on a quite
sophisticated theorem of cohomology of Lie algebras (namely, the
Hopf-Koszul-Samelson theorem, see e.g. \cite{Kostant97}) from which the \al
theorem is a consequence, modulo some combinatorial identities concerning the
trace \cite{Kostant58}. We can give a more economical proof based on similar
arguments, but that does not rely on the Hopf-Koszul-Samelson theorem. Our
proof uses only elementary properties of the Chevalley-Cartan's transgression
operator \cite{Chevalley, Cartan} and some identities concerning the
invariants $\Tr (X^k)$. It will not be presented in this paper; however, a
completely similar reasoning will allow us to handle the orthosymplectic case
$\osp (1, 2n)$.

The goal of this paper is to study possible versions of the \al
theorem in the case of Lie superalgebras. Consider the Lie super
algebra $\gl (p,q)$ and define for $X_1, \dots, X_k \in \gl (p,q)$:
\[ \Ac_k (X_1, \dots, X_k) := \sum_{\sigma \in \Sk_k} \epsilon (\sigma)
\epsilon (\sigma, \Xc) X_{\sigma(1)} \dots X_{\sigma(k)} \]
\noindent where the super sign $\epsilon(\sigma, \Xc)$ will be defined in
Section 1. The polynomial $\Ac_k$ is invariant under the action of the super
algebra $\gl (p,q)$. We call $\Ac_k$ the standard super polynomial of order
$k$ and it clear that this polynomial is a natural candidate to replace $I_k$
in the case of the superalgebra $\gl (p,q)$. The next step is to check whether
$\Ac_k$ is zero for $k$ sufficiently big. However, if $pq \neq 0$, one can
easily see that this is not true: there always exists a non nilpotent element
$X \in \gl(p,q)_{\un}$ and since $\Ac_k(X, \dots, X) = k!  X^k$, it results
that $\Ac_k \neq 0$ for all $k$. Therefore there is no standard super identity
for $\gl (p,q)$. With this counter-example in mind, one might think there is
little hope in finding such an identity for the simple subalgebras of $\gl
(p,q)$.

However, a closer look at the counter-example shows that it can be translated
in terms of invariants: the algebra of invariants of $\gl (p,q)$ is not
finitely generated \cite{Sergeev}. If we follow the philosophy of Kostant's
proofs, the algebra of invariants of the considered Lie superalgebra should
be a polynomial algebra, which leaves us with a single choice:
$\osp(1,2n)$. For this series of Lie superalgebras, the algebra of invariants
is a polynomial algebra in $n$ variables by a theorem proved by V. Kac
\cite{Musson97}. In addition, it is easy to see that all elements in $\osp (1,
2n)_\un$ are nilpotent (see Section 2), so the counter-example above does not
apply.

As a consequence, the series $\osp (1, 2n) $ seems to be a good candidate for
an \al super theorem and our goal in this paper is to show that this super
version does exists. The main result presented here is the following:

\medskip

{\bf{THEOREM:}} For $X_1, \dots, X_{4n+2} \in \osp(1, 2n)$, $\Ac_{4n+2} (X_1,
\dots, X_{4n+2}) = 0$.

\medskip

\noindent Notice that the number $4n +2$ appearing in the above theorem is
precisely the one for $\gl(2n+1)$ in the classical case of the \al theorem.

As we mentioned before, the proof of this theorem follows the lines of
Kostant's cohomological proof, but in a simpler form. Our proof does not need
to use a powerful theorem such as Hopf-Koszul-Samelson's for  $\osp
(1, 2n)$ (see \cite{Fuks-Leites}), but only elementary properties of a (super)
transgression operator and some identities concerning super traces.

We believe that in general there is no super \al theorem for the classical Lie
superalgebras, with exception made to the series $\osp (1, 2n)$. This can be
explained by the fact that their algebra of invariants is not (in general)
finitely generated. Recall that $\osp(1, 2n)$ are the only simple Lie
superalgebras (together with simple Lie algebras) that are also semi simple
\cite{Djokovic-Hochschild} (meaning complete reducibility of their
finite-dimensional representations). The fact that these Lie superalgebras
satisfy an identity of \al type strengthens the impression that they are very
close to simple Lie algebras. However, the existence of a ghost center and of
exotic primitive ideals in the enveloping algebra \cite{Musson92, Musson97,
Pinczon} indicate that the analogy cannot be carried much further.

We want to stress that the present study was performed in the context of an
invariant theory for Lie superalgebras and in that spirit. It would of course
be interesting to relate our super identity with the general theory of
PI-algebras (a very active domain, see e.g. \cite{BDDK, Raz, Rowen80}) where
the classical \al theorem plays an important role. That is a different study,
which remains to be done since our super identity does not seem to appear in
the PI-algebras literature.


\section{Notations}

\subsection{Algebras of supersymmetric and skew supersymmetric multilinear
  mappings}

Let $V = V_\ze \oplus V_\un$ be a finite-dimensional $\ZZ_2$-graded vector
space.  Considered elements $X \in V$ are supposed homogeneous, and we denote
by a small $x$ the degree. On $W = \CC$, set $W_\ze = \CC$ and $W_\un = \{ 0
\}$. Let $\Fc (V)$ be the $\ZZ$-graded space of multilinear forms on $V$ and
$\Fc^p (V)$ the subspace of $p$-forms. Consider the natural $\ZZ_2$-grading on
$\Fc^p(V)$:
\[ F \in \Fc^p (V), \deg_{\ZZ_2} (F) = f \text{ iff } \deg_{\ZZ_2} (F (X_1,
\dots, X_p)) = x_1 + \dots + x_p + f\]

The space $\Fc (V)$ is endowed with the usual tensor product $\otimes$, and
with a super tensor product denoted by $\stimes$ and defined as :
\[(F\stimes G)(X_1,\ldots,X_{p+q}) := (-1)^{g (x_1+\ldots+x_p)}
F(X_1,\ldots,X_p) G(X_{p+1}, \ldots,X_{p+q}),\]

for $X_1, \dots, X_{p+q} \in V$, $F\in\Fc^p_f(V)$, $G\in\Fc^q_g(V)$ with
$\deg_{\ZZ_2} (F) = f$ and  $\deg_{\ZZ_2} (G) = g$.

Let $\Xc = (X_1, \dots, X_p) \in V^p$ and $\sigma$ an element of the symmetric
group $\Sk_p$. Define:
\[\e(\sigma,\Xc) := (-1)^{K(\sigma,\Xc)}\]

where $K(\sigma,\Xc) := \sharp \{(i,j) \mid X_{\sigma(i)}, X_{\sigma(j)}\in
\VI, i<j \textrm{ and } \sigma(i)>\sigma(j)\}$. It follows from the definition
that $\e (\sigma, \Xc)$ is a multiplier, that is:
\[ \e(\sigma\sigma',\Xc) = \e(\sigma,\Xc) \e(\sigma',\sigma^{-1}\cdot\Xc) \] 

with $\sigma\cdot\Xc := (X_{\sigma^{-1}(1)},\ldots,X_{\sigma^{-1}(p)}).$ 

We can consider three actions of $\Sk_p$ on $\Fc^p (V)$: 
\begin{eqnarray*}
\sigma\cdot F(X_1,\ldots,X_p) &:=&
F(X_{\sigma(1)},\ldots,X_{\sigma(p)}),\label{eq1.4} \\
\sigma\sdot F(X_1,\ldots,X_p) &:=& \e(\sigma,\Xc)
F(X_{\sigma(1)},\ldots,X_{\sigma(p)}),\notag \\
\sigma\adot F(X_1,\ldots,X_p)  &:=& \e(\sigma) \e(\sigma,\Xc)
F(X_{\sigma(1)},\ldots,X_{\sigma(p)}).\label{eq1.5} 
\end{eqnarray*}

We then say that a $p$-form $F$ is {\em{supersymmetric}} if $\sigma \sdot F =
F$, $\forall \ \sigma \in \Sk_p$ and {\em{skew supersymmetric}} if $\sigma
\adot F = F$, $\forall \ \sigma \in \Sk_p$. We denote by $\Pc (V)$ the space
of supersymmetric forms and by $\Ac (V)$ the space of skew supersymmetric
forms.

Now let $\So$ and $\Ao$ be two operators on $\Fc (V)$ defined as:

\begin{equation*}\label{eq1.6}
\So(F) := \dsum_{\sigma\in\Sk_p} \sigma\sdot F, \; \; \; \; \Ao(F) :=
\dsum_{\sigma\in\Sk_p} \sigma\adot F, \; \; \; \; \forall \ F\in \Fc^p(V).
\end{equation*}

We can then define a product on $\Pc (V)$ and $\Ac (V)$ as:
\begin{equation*}\label{eq1.7}
F \cdot G := \dfrac1{p!q!} \So(F\stimes G),
\end{equation*}
for $F\in\Pc^p(V)$, $G \in \Pc^q(V)$,
\begin{equation*}\label{eq1.8}
F\wedge G := \dfrac1{p!q!} \Ao(F\stimes G), 
\end{equation*}
for $F\in\Ac^p(V)$, $G \in \Ac^q(V)$.  

This gives an algebra structure on $\Pc(V)$ and $\Ac(V)$. The algebra $\Pc
(V)$ is $\ZZ_2$-graded (since $V$ is $\ZZ_2$-graded) and isomorphic to the
(usual) tensor product $\sym(\VO^*) \otimes \ext (\VI^*)$. The algebra $\Ac
(V)$ is double graded by $\ZZ \times \ZZ_2$, and isomorphic to $\displaystyle
\ext (\VO^*) \! \! \underset{\ZZ \times \ZZ_2}{\otimes} \! \! \sym(\VI^*)$. We
have
\[F\cdot G = (-1)^{fg}G \cdot F\]
for $F,G\in\Pc(V)$, $\deg_{\ZZ_2} (F) = f$, $\deg_{\ZZ_2}(G) = g$, and
\[F\wedge G = (-1)^{nm+fg} G\wedge F,\] 
for $F,G \in\Ac(V)$, $\deg_{\ZZ\times\ZZ_2}(F) = (n,f)$,
$\deg_{\ZZ\times\ZZ_2}(G) = (m,g)$.

These relations imply that $\Pc(V)$ and $\Ac (V)$ are supercommutative with
respect to their gradation. We can say that $\Pc (V)$ (respectively $\Ac (V)$)
is the analogous of the algebra of polynomial functions (respectively of the
Grassman algebra) in the non graded case.

The following formulae will be useful in this work: let $\phi_1, \dots, \phi_p
\in V^*$ with $\ZZ_2$-degrees $\varphi_1, \dots, \varphi_p$, $\varphi := (
\varphi_1, \dots, \varphi_p)$, then
\begin{eqnarray*}
&\phi_1 \cdot \ldots \cdot \phi_p = (-1)^{\Omega(\varphi,\varphi)} \So
  (\phi_1\otimes \ldots\otimes\phi_p) \  \text{     and } \\
&\phi_1\wedge \ldots \wedge \phi_p = \Ao (\phi_1\stimes\ldots\stimes\phi_p) =
  (-1)^{\Omega(\varphi,\varphi)} \Ao (\phi_1\otimes\ldots\otimes\phi_p), 
\end{eqnarray*}

where $\Omega$ is the $2$-form with matrix $\begin{pmatrix} 0 & \ldots &
  \ldots & 0 \\ 1& \ddots & 0 & \vdots \\ \vdots & \ddots & \ddots & \vdots \\
  1 & \ldots & 1 & 0 \\ \end{pmatrix}$.

For $X \in V$, define super derivations $\Do_X$ and $\io_X$ of $\Pc(V)$ and
$\Ac (V)$ respectively as:
\[\Do_X(F)(X_1,\ldots,X_{p-1}) := (-1)^{xf} F(X,X_1, \ldots, X_{p-1})\]
for $F\in\Pc (V)$, $\deg_{\ZZ_2}(F) = f$ and
\[\io_X(F)(X_1,\ldots,X_{p-1}) := (-1)^{xf} F(X,X_1, \ldots,X_{p-1}) \]
for $F\in\Ac(V)$, $\deg_{\ZZ\times\ZZ_2}(F) = (p,f)$.

Hence, $\Do_X$ is a super derivation of degree $x$ of $\Pc (V)$ and $\io_X$ is
a super derivation of degree $(-1, x)$ of $\Ac(V)$.

\subsection{Cohomology of Lie superalgebras (see \cite{Leites})}

Let $\g = \gO \oplus \gI$ be a Lie superalgebra with $\dim \gO = p$ and $\dim
\gI = q$. The contragredient representation $\check{\ad}$ of the adjoint
representation $\ad$ can be extended to a representation $\Lo^s$ of $\g$ into
$\Pc(V)$ and to a representation $\Lo^a$ of $\g$ into $\Ac(V)$. For $X \in
\g$, $\Lo_X^s$ (resp. $\Lo_X^a$) is the super derivation of degree $x$ (resp.
$(0,x)$) of $\Pc(V)$ (resp. $\Ac (V)$) defined as: for $F \in \Pc(\g)$ (resp.
$\Ac(\g)$) with $\deg_{\ZZ_2}(F) = f$ (resp. $\deg_{\ZZ\times\ZZ_2}(F) =
(n,f)$),

\[ \Lo_X^{a,s}F(X_1,\ldots,X_n) := -(-1)^{xf} \dsum_{j=1}^n
(-1)^{x(x_1+\ldots+x_{j-1})} F(X_1,\ldots,\ad X(X_j),\ldots,X_n).\]

Denote by $I^s ( \g)$ and $I^a (\g)$ the invariants under these actions. Let
$d$ be the map from $V^*$ to $\Ac(V)$ defined as:

\begin{equation*}\label{eq1.16}
d\phi(X_1,X_2) := -\phi([X_1,X_2]), \forall \ \phi\in \g^*.
\end{equation*}

There exists a super derivation (also denoted by $d$) of $\Ac(V)$ of degree
$(1,0)$ extending $d$: for $F\in \Ac(\g)$,

\begin{eqnarray*}\label{eq1.17}
dF(X_1,\ldots,X_{n+1}) := \dsum_{i<j} (-1)^{i+j} &(-1)^{x_i(x_1+\ldots+
  x_{i-1})} (-1)^{x_j(x_1+\ldots+\widehat{x_i}+\ldots+x_{j-1})}\notag{} \\ &
  F([X_i,X_j],X_1,\ldots,\widehat {X_i},\ldots,\widehat{X_j},\ldots,X_{n+1}).
\end{eqnarray*}

{}From the Jacobi identity, it comes $d^2 = 0$ and we can then define the
cohomology (with trivial coefficients) of $\g$ as:
\[ Z(\g) := \mathrm{Ker}(d), \ B(\g) := \mathrm{Im}(d) \textrm{ and } H(\g) :=
Z(\g)/B(\g).\]

Let $\{ X_1, \dots, X_{p+q} \}$ be a basis of $\g$ and $\{ \phi_1, \dots,
\phi_{p+q} \}$ its dual basis. Define the forms $\tilde{\phi}_i$ as
$\tilde{\phi}_i (X) := (-1)^{x_i x} \phi_i (X)$, $X \in \g$. Thus, one has:

\begin{equation}\label{eq1.18}
d = \dfrac12 \dsum_{i=1}^{p+q} \tilde{\phi_i}\wedge \Lo_{X_i}^a.
\end{equation}

It results from (\ref{eq1.18}) that $I^a (\g) \subset Z(\g)$. Moreover, one
has: 
\begin{equation}\label{eq1.19}
\Lo_X^a = \io_X\circ \ d + d\circ \io_X, \ \forall \ X\in\g.
\end{equation}

As a consequence, $\Lo_X^a$ commutes with $d$ and $\Lo_X^a (Z (\g)) \subset
B(\g)$.


\section{Orthosymplectic Lie superalgebras}

In this section, let $\g$ be the orthosymplectic Lie superalgebra $\osp(1,
2n)$.   Among simple Lie superalgebras, the orthosymplectic
$\osp(1,2n)$ are the only ones (together with simple Lie algebras) satisfying
the remarkable property of being semi simple \cite{Djokovic-Hochschild} ,
meaning that every finite-dimensional representation is completely reducible.

\subsection{The Weyl algebra and $\osp(1,2n)$} 

In the quantization framework, the Lie superalgebra $\g$ can be realized as
follows: let $\Ab_n$ be the Weyl algebra generated by $\{ p_i, q_i, i = 1,
\dots, n \}$ with $[p_i,q_i]_\Lc = 1$, $\forall \ i$, $[p_i,q_j]_\Lc =
[p_i,p_j]_\Lc = [q_i,q_j]_\Lc = 0$, if $i\ne j$ where $[ \cdot, \cdot]_\Lc $
denotes the Lie bracket. The algebra $\Ab_n$ is $\ZZ_2$-graded, hence a Lie
superalgebra. Denote by $[ \cdot, \cdot]$ its bracket.

\begin{defn} The twisted adjoint action of $\Ab_n$ onto itself is defined as:
\begin{equation*}\label{eq2.1}
\ad'A(B) := AB-(-1)^{a(b+1)} BA
\end{equation*}

for $A,B\in\Ab_n$, $\deg_{\ZZ_2}(A) = a$, $\deg_{\ZZ_2}(B) = b.$

\end{defn}

Let $\VI := \spa \{ p_i, q_i, i = 1, \dots, n \}$ and $\hk := \VI \oplus [\VI,
\VI]$. Then $\hk$ is a subalgebra of the Lie superalgebra $\Ab_n$. Let now $V
:= \VO \oplus \VI$ where $\VO := \CC \cdot 1$. We have $\ad' \hk (V) \subset
V$.  Moreover the supersymmetric $2$-form $F (X,Y) := [X, Y]_\Lc$, $X, Y \in
\VI$ and $F(1,1) := -2$ is $\ad' \hk$-invariant. It follows that $\hk \simeq
\osp(1, 2n)$. An easy but remarkable consequence is the following

\begin{prop}\label{pr2.1}
If $X \in \osp(1, 2n)_\un$, then $X^3= 0$.
\end{prop}

\begin{proof}
  It is enough to show that if $X \in \VI$, then $(\ad' X |_V)^3 = 0$. Using
   $(\ad' X )(1) = 2 X$ and $(\ad' X)^2 (Y) = 2 [X,Y]_\Lc \ X$,
  $\forall \ Y \in \VI$, the result follows.
\end{proof}

More generally:

\begin{prop}
Let $\pi$ be a finite-dimensional representation of  $\g = \osp(1,
2n)$. If $X \in \gI$, then $\pi (X)$ is nilpotent.
\end{prop}

\begin{proof}
  We use here the realization of $\g$ as $\hk$. Let $X \in \hk_\un = \VI$, $X
  \neq 0$. There exists a Darboux basis of $\VI$ for the form $F|_{\VI \times
    \VI}$ such that $X$ is the first basis element. We can then suppose that
  $X = p_1$. Let $\lk = \lk_\ze \oplus \lk_\un$ with $\lk_\un = \spa \{ p_1,
  q_1 \}$ and $\lk_\ze = [ \lk_\un, \lk_\un ]$. So $\lk \simeq \osp(1,2)$. Let
  $\rho = \pi|_\lk$.  Write $\rho = \oplus_{i \in I} \ \rho_i$ its
  decomposition into simple components. If $d_0 = \max \{ \dim \rho_i, i \in I
  \}$, then $\pi (p_1)^{d_0} = 0$.
\end{proof}

\subsection{Cohomology of $\osp(1,2n)$}

{}From \cite{Djokovic-Hochschild}, the representation $\Lo^a$ of $\g$ is
completely reducible. Using Koszul' strategy in \cite{Koszul}, this fact
together with the results in Section 1.2, in particular the equations
(\ref{eq1.18}) and (\ref{eq1.19}), allows us to prove

\begin{lem}
Every cohomology class of $H (\g)$ contains one and only one invariant
cocycle. In particular, $0$ is the unique invariant coboundary and $H(\g) =
I^a (\g)$.
\end{lem}

For the sake of completeness, we should mention that there exist better
results concerning $H(\g)$: Fuks and Leites \cite{Fuks-Leites} have announced
that $H(\g) \simeq H(\g_\ze) = H (\spk(2n))$. However, we shall not need these
results here.

\subsection{Invariants}

Concerning $I^s (\g)$, it results from V. Kac's work that the Chevalley
  restriction theorem holds \cite{Musson97}: let $\hk$ be a Cartan subalgebra
  of $\g_\ze$ and $W$ the Weyl group, then the restriction of $I^s (\g)$ into
  $\sym(\hk^*)^W$ is an algebra isomorphism. As a consequence, $I^s (\g)$ is a
  polynomial algebra in $n$ variables. We will see later how to choose
  convenient generators.


\section{Chevalley's transgression operator for Lie superalgebras}

The transgression operator $t \colon \sym(\g^*) \rightarrow \ext (\g^*)$ was
introduced by Chevalley (\cite{Chevalley, Cartan}, see also \cite{Kostant97})
and it is a fundamental tool in the theory of Lie algebras. In this section,
we shall generalize this notion to the case of Lie superalgebras and give some
elementary properties that will be useful in the sequel.

Let $\g = \g_\ze \oplus \g_\un$ be a Lie superalgebra. Let $\{ X_1, \dots,
X_p\}$ be a basis of $\g_\ze$, $\{ Y_1, \dots, Y_q \}$ a basis of $\g_\un$,
$\{\Omega_1, \dots, \Omega_p\}$ and $\{ \phi_1, \dots, \phi_q \}$ their
respective dual basis. There exists a super derivation $\Ro$ of $\Pc (\g)$ of
degree $0$ extending $\Id_{\g^*}$:
\[ \Ro := \dsum_{i=1}^p \Omega_i \Do_{X_i} - \dsum_{j=1}^q \phi_j \Do_{Y_j} \]

We have:
\begin{equation*}\label{3.1}
\Ro(P) = (\deg_\ZZ P) \ P, \quad \forall \ P\in \Pc(\g)
\end{equation*}

where $\deg_\ZZ P$ comes from $\Pc (\g) = \sym (\gO^*) \otimes \ext (\gI^*)$
and from the natural $\ZZ$-gradations of $\sym( \gO^*)$ and $\ext (\gI^*)$.

There exists an algebra homomorphism $s \colon \Pc (\g) \rightarrow \Ac (\g)$
such that $s (\Omega_i) = d \Omega_i$, $i = 1, \dots, p$, and $s (\phi_j) = d
\phi_j$, $j = 1, \dots, q$ (since the $d \Omega_i$ ($i = 1, \dots, p$)
commute, the $d \phi_j$ ($j = 1, \dots, q$) anticommute and the $d \Omega_i$,
$d \phi_j$ ($i= 1, \dots, p$, $j = 1, \dots, q$) commute).

One can easily check that $d (s(P))=0$, $\forall \ P \in \Pc (\g)$. Besides,
$s$ is a homomorphism of $\g$-modules, if $\Pc (\g)$ is endowed with the
representation $\Lo^s$ and $\Ac (\g)$ endowed with the representation
$\Lo^a$. Therefore $s(I^s (\g)) \subset I^a (\g)$.

Following Chevalley, we now set:

\begin{defn}
The transgression operator $t \colon \Pc (\g) \rightarrow \Ac (\g)$  is
defined as
\begin{equation}\label{eq3.2}
t(P) := \dsum_{i=1}^p \Omega_i\wedge s(\Do_{X_i}(P)) - \dsum_{j=1}^q
\phi_j\wedge s(\Do_{Y_j}(P)), \forall \ P\in\Pc(\g)
\end{equation}
\end{defn}

A priori, this definition seems to be basis dependent, but this is not the
case as we shall show below. For the time, let us state:

\begin{lem}\label{lem3.1}
One has $d(t(P)) = s(\Ro(P))$, $\forall \ P \in \Pc (\g)$.
\end{lem}

Since $\Ro (P) = (\deg_\ZZ P) \ P$, Lemma \ref{lem3.1} shows that if $P$ has
no constant term, then $s(P)$ is a coboundary.

Moreover $t$ is an $s$-derivation:

\begin{lem}\label{lemtder}
One has $t(P \cdot Q) = t(P)\wedge s(Q)+s(P)\wedge t(Q)$, for all
$P,Q\in\Pc(\g)$.
\end{lem}

In order to establish some other properties of the transgression, we need now
an intrinsic definition of $t$. First, observe that there is an isomorphism
$\End (\g) = \g^* \stimes \g$ given by:

\begin{equation*}\label{eq3.4}
(\Omega\stimes X)(Y) := (-1)^{xy} \Omega(Y) \ X, \forall \ \Omega\in\g^*, \
  X,Y\in\g.  
\end{equation*}

Thanks to this identification, the representation $\pi := \check{\ad} \otimes
\ad$ becomes $\ad (\ad \cdot)$ and $\Id_\g = \dsum_{i=1}^p \Omega_i\stimes X_i
- \dsum_{j=1}^q \phi_j\stimes Y_j$ is $\pi$-invariant.

Now fix $P \in \Pc (\g)$ and set $\tau_p \colon \End (\g) \rightarrow \Ac (\g)$
as

\begin{equation}\label{eq3.5}
\tau_P(\Omega\stimes X) := \Omega\wedge s(\Do_X(P))
\end{equation}

It is immediate that $\tau_p (\Id_\g) = t(P)$, so the definition of $t$ in
(\ref{eq3.2}) is basis independent. In addition, using the representation
$\pi$ on $\End (\g)$ and $\Lo^a$ on $\Ac(\g)$, one has

\begin{lem}\label{lem3.2}
If $P \in I^s (\g)$, then $\tau_P \colon \End (\g) \rightarrow \Ac (\g)$ is a
$\g$-module homomorphism.
\end{lem}

As a direct consequence of (\ref{eq3.5}) and Lemma \ref{lem3.2}, we obtain:
\begin{equation}\label{eq3.6}
t(I^s(\g)) \subset I^a(\g)
\end{equation}

Combining \ref{eq3.6}, \ref{eq1.18} and Lemma \ref{lem3.1}, one has

\begin{lem}\label{lem3.3}
Let $I^s_+ (\g)$ be the subspace of $I^s (\g)$ with no constant terms. Then
for all $P \in I^s_+ (\g)$, $s(P) = 0$.
\end{lem}

Finally applying Lemma \ref{lemtder}, we conclude

\begin{lem}\label{lem3.4}
For all $P\in \CC \oplus (I_+^s(\g))^2$, $t(P) = 0$.
\end{lem}

\begin{rem}
  For similar results in the non graded case, see \cite{Chevalley} or
  \cite{Kostant97}.
\end{rem}


\section{Standard super polynomials and super identities in $\gl (p,q)$}

In this section, $V = \VO \oplus \VI$ with $\dim \VO = p$, $\dim \VI = q$, and
$\g$ is the Lie superalgebra $\g = \End (V) \simeq \gl (p,q)$.

We identify $\End (V)$ and $V \otimes V^*$ by using:
\begin{equation*}\label{eq4.1}
Z\otimes \Omega(T) := Z.\Omega(T), \ \forall \ Z,T\in V, \ \Omega\in V^*
\end{equation*}

Then define the super trace on $\g$ as:
\begin{equation*}\label{eq4.2}
\str(Z\otimes\Omega) := (-1)^{\omega z}\Omega(Z), \ \forall \ Z\in V, \
\Omega\in V^*
\end{equation*}

\begin{rem}\label{rem4.1}
  With this definition, the $2$-form $B (Z | T) := \str(ZT)$ is supersymmetric
  and non degenerate on $\g$. In the case $p = 1$ and $q = 2n$,
  $B|_{\osp(1,2n)}$ is non degenerate as well.
\end{rem}

\begin{defn} 
The standard supersymmetric super polynomials $\Pc_k$ (resp. skew
supersymmetric $\Ac_k$) are given by:
\begin{subequations}
\begin{eqnarray*}
\Pc_k(X_1,\ldots,X_k) &:=& \dsum_{\sigma\in\Sk_k} \e(\sigma;\Xc)
X_{\sigma(1)}\ldots X_{\sigma(k)},\label{eq4.3a} \\  
\Ac_k(X_1,\ldots,X_k) &:=& \dsum_{\sigma\in\Sk_k} \e(\sigma) \e(\sigma;\Xc)
X_{\sigma(1)}\ldots X_{\sigma(k)},\label{eq4.3b}
\end{eqnarray*}
\end{subequations}
where $k\ge 1$, $X_1, \dots, X_k \in \g$.

\end{defn}

The polynomials $\Pc_k$ and $\Ac_k$ are $\g$-invariant $k$-linear maps from
$\g^k$ to $\g$. They verify the recursive relations below:
\begin{subequations}
\begin{equation}\label{eq4.4a}
\Pc_{k+1}(X_1,\ldots,X_{k+1}) = \dsum_{j=1}^{k+1}
(-1)^{x_j(x_1+\ldots+x_{j-1})} X_j. 
\Pc_k(X_1,\ldots,\hat{X_j},\ldots,X_{k+1}),
\end{equation}
\begin{eqnarray}\label{eq4.4b}
 \Ac_{k+1}(X_1,\ldots,X_{k+1}) &=& \dsum_{j=1}^{k+1} (-1)^{j+1}
 (-1)^{x_j(x_1+\ldots+x_{j-1})} X_j.\notag\\ 
& & \Ac_k(X_1,\ldots,\hat{X_j},\ldots,X_{k+1}).
\end{eqnarray}
\end{subequations}

{}From $\Pc_k$ and $\Ac_k$, we can construct $P_k \in I^s(\g)$ and $\Lambda_k
\in I^a (\g)$:
\begin{subequations}
\begin{eqnarray*}
P_k(X_1,\ldots,X_k) &:=& \str(\Pc_k(X_1,\ldots,X_k)),\label{eq423a} \\
 \Lambda_k(X_1,\ldots,X_k) &:=& \str(\Ac_k(X_1,\ldots,X_k))\label{eq423b}
\end{eqnarray*}
\end{subequations}

\begin{prop}\label{pr4.1}
One has:\\
\noindent {\it{(a)}} \begin{eqnarray}
& & P_{2k+1}(X_1,\ldots,X_{2k+1}) = (2k+1)
  B(\Pc_{2k}(X_1,\ldots,X_{2k})|X_{2k+1}),\notag \\
& & \Lambda_{2k}(X_1,\ldots,X_{2k}) = 0, \label{eq4.6} \\
& & \Lambda_{2k+1}(X_1,\ldots,X_{2k+1}) = (2k+1)
  B(\Ac_{2k}(X_1,\ldots,X_{2k})|X_{2k+1}).\notag 
\end{eqnarray}
\noindent{\it{(b)}}\begin{equation}\label{eq4.7}
\dsum_{\sigma\in\Sk_{2k}}\e(\sigma) \e(\sigma;\Xc)
     [X_{\sigma(1)},X_{\sigma(2)}]\ldots [X_{\sigma(2k-1)},X_{\sigma(2k)}] =
     2^k \Ac_{2k}(X_1,\ldots,X_{2k}) 
\end{equation}
\noindent{\it{(c)}}\begin{eqnarray}
& & \dsum_{\sigma\in\Sk_{2k+1}} \e(\sigma) \e(\sigma;\Xc)
[X_{\sigma(1)},X_{\sigma(2)}]\ldots [X_{\sigma(2j-1)},X_{\sigma(2j)}]
X_{\sigma(2j+1)} \notag\\
& & [X_{\sigma(2j+2)},X_{\sigma(2j+3)}] \ldots
[X_{\sigma(2k)},X_{\sigma(2k+1)}] = 2^k \Ac_{2k+1}(X_1,\ldots,X_{2k+1})
\notag\\ 
\label{eq4.8}  
\end{eqnarray}
\end{prop}

\begin{rem} 
The identities (\ref{eq4.6}), (\ref{eq4.7}) and (\ref{eq4.8}) are super
  versions of classical identities in the non graded case. Their proofs are
  simple adaptations to the super case. Other super identities can be settled,
  but they will not be needed in this work.
\end{rem}

Let us examine what happens when we apply the transgression on the invariant
$P_k$ defined by the super trace.

\begin{thm}\label{thm4.1}
  One has $t(P_k) = (-1)^{k-1}k \Lambda_{2k-1}$.
\end{thm}

\begin{proof} The main argument here will be Lemma \ref{lemtder}.  Let
  $M_{ij}$ be the coordinate forms. Then 
\[M_{ii}(X_1\ldots X_k) =
\dsum_{R=(r_1,\ldots,r_{k-1})} (-1)^{\Omega(m_{iR},m_{iR})}
M_{ir_1}\stimes\ldots\stimes M_{r_{k-1}i} (X_1,\ldots,X_k)\] 
where $m_{iR} :=\begin{pmatrix} m_{ir_1} \\ \vdots \\ m_{r_{k-1}i}
\end{pmatrix}$.

Supersymmetrizing, we obtain:
\begin{eqnarray*}
P_k &=& \dsum_{i\in\dc{1,p} \atop R} (-1)^{\Omega(m_{iR},m_{iR})} M_{ir_1}
\cdot M_{r_1r_2}\cdot \ldots \cdot M_{r_{k-1}i} \\ & & -
\dsum_{j\in\dc{p+1,p+q} \atop R} (-1)^{\Omega(m_{jR},m_{jR})} M_{jr_1} \cdot
M_{r_1r_2} \cdot \ldots \cdot M_{r_{k-1}j}
\end{eqnarray*}

\noindent (notice that the products above are calculated in $\Pc(\g)$).

{}From $t(M_{rs}) = M_{rs}$, $\forall \ r,s$ and Lemma \ref{lemtder}, it
comes:
\[ t(M_{ir_1} \cdot \ldots \cdot M_{r_{k-1}i}) = \dsum_{\ell=1}^k
dM_{ir_1}\wedge dM_{r_1r_2}\wedge \ldots\wedge
M_{r_{\ell-1}r_\ell}\wedge\ldots \wedge dM_{r_{k-1}i}\] 
(if $\ell=k$ then $r_k = i$ in the sum).

Therefore:
\begin{eqnarray*}
& & t(M_{ir_1} \cdot \ldots \cdot M_{r_{k-1}i})(X_1,\ldots,X_{2k-1}) \\ 
&=& (-1)^{\Omega(m_{iR},m_{iR})} \dfrac{(-1)^{k-1}}{2^{k-1}}
  \dsum_{\sigma,\ell} \e(\sigma) \e(\sigma,\Xc)
  M_{ir_1}([X_{\sigma(1)},X_{\sigma(2)}]) \ldots \\ 
& & M_{r_{\ell-1}r_\ell}(X_{\sigma(2\ell-1)})\ldots
  M_{r_{k-1}i}([X_{\sigma(2k-2)},X_{\sigma(2k-1)}]) 
\end{eqnarray*}

At the end, we have:
\begin{eqnarray*}
& & \dsum_{R} (-1)^{\Omega(m_{iR},m_{iR})} t(M_{ir_1} \cdot \dots \cdot
  M_{r_{k-1}i})(X_1,\ldots,X_{2k-1}) \\ 
&=& \dfrac{(-1)^{k-1}}{2^{k-1}} \dsum_{\sigma,R,\ell} \e(\sigma)
  \e(\sigma,\Xc) M_{ir_1}([X_{\sigma(1)},X_{\sigma(2)}]) \ldots \\
& & M_{r_{\ell-1}r_\ell}(X_{\sigma(2\ell-1)})\ldots
  M_{r_{k-1}i}([X_{\sigma(2k-2)},X_{\sigma(2k-1)}]) \\ 
&=& (-1)^{k-1} \dsum_\ell M_{ii} (\Ac_{2k-1}(X_1,\ldots,X_{2k-1})) \qquad
\text{ (by \ref{eq4.8})}\\
&=& (-1)^{k-1} k M_{ii} (\Ac_{2k-1}(X_1,\ldots,X_{2k-1})).
\end{eqnarray*}

\end{proof}


\section{The Amitsur Levitzki theorem for $\osp(1,2n)$}

Henceforth we will assume that $\g = \osp(1,2n)$ and $\gt = \gl(1,2n)$. We
will now prove a (super) version of the \al theorem for $\g$. In other words,
we will show:

\begin{thm}\label{thm5.1}
For all $X_1,\ldots,X_{4n+2}\in\g$, we have $\Ac_{4n+2}(X_1,\ldots,X_{4n+2}) =
0.$
\end{thm}

Notice that this identity is valid if $X_1,\dots, X_{4n+2} \in \gO$ by the
classical  \al theorem. Furthermore, if $X_1 = \dots = X_{4n+2} = X
\in \gI$ then by Proposition \ref{pr2.1}, the identity holds as well.

The theorem will be a consequence of Theorem \ref{thm4.1} and two lemmas:

\begin{lem}\label{lem5.1}

One has: 

\begin{itemize}

\item[(1)] For all $X_1,\ldots,X_{2p+1}\in\g$, $\Pc_{2p+1}(X_1,\ldots,X_{2p+1})
  \in \g$.

\item[(2)] For all $X_1,\ldots,X_{4p+1}\in\g$,
  $\Ac_{4p+1}(X_1,\ldots,X_{4p+1}) \in \g$.
  
\item[(3)] For all $X_1, \dots, X_{4p+2} \in \g$, $\Ac_{4p+2} (X_1, \dots,
  X_{4p+2}) \in \g$.

\end{itemize}

\end{lem}

As a consequence, $P_{2k+1}$, $\Lambda_{4p+1}$ and $\Lambda_{4p+2}$ vanish as
multilinear mappings on $\g$.

Recall from Subsection 2.3 that the restriction $R \colon I^s (\g) \rightarrow
J$ is an algebra isomorphism where $J := \sym(\hk^*)^W$. The elements of $\hk$
are the matrices  $H(\alpha_1,\alpha_2,\ldots,\alpha_n) =$ $
\begin{pmatrix} \alpha_1 & 0 & & & & \\ 0 &-\alpha_1 & & & 0 & \\ & & \ddots &
  & & \\ & 0 & & & \alpha_n & 0 \\ & & & & 0 &-\alpha_n \\ \end{pmatrix}$, so
one has $\sym(\hk^*) =$  $\CC[\alpha_1, \dots, \alpha_n]$ and the
Weyl group is generated by permutations and changes signs of $\alpha_1, \dots,
\alpha_n$. For these reasons, we can write $J = \CC [t_1, \dots, t_n]$ where
$t_k := \sum_{i =1}^k \alpha_i^{2k}$ for $1 \leq k \leq n$. It is clear that
$t_k \in J$, $\forall \ k$ and that $t_k \in J_+^2$ if $k \geq n+1$ where
$J_+$ denotes the augmentation ideal. On the other hand, $R (P_{2k}) = 2 t_k$,
therefore one deduces:

\begin{lem}
One has $I^s(\g) = \CC[P_2,P_4,\ldots,P_{2n}]$ and $P_{2n+2}\in
(I_+^S(\g))^2$. 
\end{lem}

We will next terminate the proof of Theorem \ref{thm5.1}.

\begin{proof} {\it(of Theorem \ref{thm5.1})}
  Let $t_\g$ be the transgression defined on $\g$ and $t_{\gt}$ be
  transgression defined on $\gt$.  Since $\g$ is a subalgebra of $\gt$, if $P$
  is a $p$-form in $\Pc(\gt)$, one has $t_{\gt} (P)|_{\g^p} = t_\g
  (P|_{\g^p})$. In the sequel, we use $t$ for both transgressions $t_\g$ and
  $t_{\gt}$, and we consider multilinear mappings restricted to $\g$. Now,
  since $P_{2n+2} \in (I^s_+ (\g))^2$, we have $t(P_{2n+2}) = 0$ from Lemma
  \ref{lem3.4}. Using Theorem \ref{thm4.1}, we deduce $t(P_{2n+2}) = -(2n+2)
  \Lambda_{4n+3}$, hence $\Lambda_{4n+3} =0$. From Proposition \ref{pr4.1},
  for all $X_1,\dots, X_{4n+3} \in \g$,
\[\Lambda_{4n+3}(X_1,\ldots,X_{4n+3}) = (4n+3)
B(\Ac_{4n+2}(X_1,\ldots,X_{4n+2})|X_{4n+3}).\]

But $\Ac_{4n+2} (X_1, \dots, X_{4n+2}) \in \g$ by Lemma \ref{lem5.1} {\em(3)},
hence from Remark \ref{rem4.1}:
\[\Ac_{4n+2}(X_1,\ldots,X_{4n+2}) = 0, \text{ for all }
X_1,\ldots,X_{4n+2}\in\g.\]
\end{proof}

\begin{rem}
{} From (\ref{eq4.4b}), we have ${\Ac_k|}_{\g^k} = 0$ if $k \geq 4n+2$. Also
  one can check that ${\Ac_{4n}|}_{\scriptstyle{\gO^{{4n-1}}\times \gI}} \neq
  0$ (thanks to Hopf-Koszul-Samelson theorem for $\gO = \spk(2n)$). So the
  index obtained in Theorem \ref{thm4.1} is the best possible, if one
  considers only even indices, a technical but justified assumption (see
  \cite{Kostant81}).  As for ${\Ac_{4n+1}|}_{\g^{4n+1}}$, it does not vanish
  if $n=1$ et $n=2$, but the general case is still to be done.
\end{rem}

\section{Acknowledgements} 
We wish to thank D.~Sternheimer and the referees for suggestion of several
improvements to this paper.

\

This article is dedicated to Moshe Flato (1937--1998). He introduced us to the
Amitsur-Levitzki theorem (among many other subjects). When Moshe was a first
year student at the Hebrew University in Jerusalem, Jakob Levitzki taught him
algebra and (after his sudden death in 1956) was replaced in the middle of the
year by his assistant Shimshon Avraham Amitsur, with whom Moshe considered for
a while to work for a PhD before choosing physics with Giulio Racah. Moshe was
an extraordinary mathematician and physicist, with insight in both fields. But
above all he was a superb human being and a great friend.

\end{document}